\documentclass[12pt]{amsart}
\usepackage{amsfonts, amssymb, amscd, amsthm, amsmath, xypic}
\usepackage{geometry}
\usepackage{graphics}
\renewcommand{\phi}{\varphi}

\begin{document}
\sloppy

\title{A class of Toeplitz operators \\ with hypercyclic subspaces}
\author{Andrei Lishanskii}

\address{ 
 Andrei Lishanskii, 
\newline
Department of Mathematics and Mechanics,
St. Petersburg State University, 
\newline
\phantom{x}\,\, and
\newline 
Chebyshev Laboratory,
St. Petersburg State University,
\newline 
St. Petersburg, Russia,
\newline {\tt lishanskiyaa@gmail.com}
\newline\newline \phantom{x}
}
\thanks{Author were supported by the Chebyshev Laboratory (St.Petersburg State University) under RF Government grant 11.G34.31.0026 and by JSC "Gazprom Neft".}

\begin{abstract}

We use a theorem by Gonzalez, Leon-Saavedra and Montes-Rodriguez to construct 
a class of coanalytic Toeplitz operators 
which have an infinite-dimensional closed subspace, where 
any non-zero vector is hypercyclic.

\end{abstract}
\maketitle

\section{Introduction}

Let $X$ be a separable Banach space (or a Frechet space), 
and let $T$ be a bounded linear operator in $X$. 
If there exists $x \in X$ such that the set 
$\{T^n x, n\in\mathbb{N}_0\}$ is dense in $X$, 
then $T$ is said 
to be a \textit{hypercyclic operator} and $x$ is called its \textit{hypercyclic vector}.
Here $\mathbb{N}_0 = \mathbb{N}\cup\{0\}$.

Dynamics of linear operators and, as a special case, theory of hypercyclic operators
was actively developed for the last 20 years. A detailed review of the 
results up to the end of 1990-s 
is given in the paper \cite{Gro}. For a recent exposition of the theory 
see the monographs \cite{Bay-Mat, GEP}.

However, first examples of hypercyclic operators appeared much earlier. 
In 1929 Birkhoff has shown that the translation operator
$T_a : f(z) \mapsto f(z + a), a \in \mathbb{C}$, 
$a \neq 0$, is hypercyclic in the Frechet space 
of all entire functions $Hol (\mathbb{C})$ with topology of uniform convergence
on the compact sets. 
Later, McLane proved hypercyclicity of the differentiation operator $D : f \mapsto f'$ on $Hol (\mathbb{C})$.
The first example of a hypercyclic operator in the Banach setting was given 
in 1969 by Rolewicz \cite{rol} who showed that for any $\lambda \in \mathbb{C}$, 
$|\lambda| > 1$, the operator $\lambda S^*$ is hypercyclic on 
$\ell^p(\mathbb{N}_0), 1 \leq p < \infty$, 
where $S^*$ is the backward shift on $\ell^p (\mathbb{N}_0)$ 
which transforms a vector $x = (x_0, x_1, \ldots, x_n, \ldots) 
\in \ell^p (\mathbb{N}_0)$ to the vector $(x_1, x_2, \ldots, x_{n+1}, \ldots)$.
\medskip

Given a hypercyclic operator $T$, what can be said about the set of 
its hypercyclic vectors? Clearly, if $x$ is a hypercyclic vector for 
the operator $T$ then $Tx, T^2 x, T^3 x, \ldots$ are hypercyclic 
vectors for $T$ as well. Hence, 
the set of hypercyclic vectors is dense when it is non-empty.

The following result was proved by Bourdon \cite{Bou}
(a special class of operators commuting with generalized backward shifts 
was previously considered by Godefroy and Shapiro in \cite{God-Sha}).
\bigskip
\\
\textbf{Theorem (Bourdon, \cite{Bou}).}
\textit{Let $T$ be a hypercyclic operator acting on a Hilbert space $H$. 
Then there exists a dense linear subspace, where any non-zero vector 
is hypercyclic for $T$.}
\bigskip
\\
\textbf{Definition.} Given a hypercyclic operator $T$, 
an infinite-dimensional closed subspace in which every non-zero vector 
is hypercyclic for $T$ is called a \textit{hypercyclic subspace}.
\bigskip

Montes-Rodriguez \cite[Theorem 3.4]{MR}
proved that the operator $\lambda S^*$, $|\lambda| > 1$, 
on $\ell^2 (\mathbb{N}_0)$ has no hypercyclic subspaces.
However, for some class of functions of the backward shift $S^*$ on 
$\ell^2 (\mathbb{N})$ there exists a hypercyclic subspace, 
and it is the main result of the present paper. 
To state it, we need to introduce some notations.
Let $\mathbb{D} = \{z\in\mathbb{C}: |z|<1\}$ 
be the unit disc and let
$\mathbb{T} = \{z\in\mathbb{C}: |z|=1\}$ be the unit circle. 
Recall that the {\it disc algebra} $A(\mathbb{D})$ 
is the space of all functions which are continuous in the closed disc
$\overline{\mathbb{D}}$ and analytic in $\mathbb{D}$ (with the 
norm $\max_{z\in \overline{\mathbb{D}}} |\phi(z)|$).
\bigskip
\\
\textbf{Main Theorem.} 
\textit{For any function $\phi \in A(\mathbb{D})$ 
such that $\phi(\mathbb{T}) \cap \mathbb{T} \ne \emptyset$ and 
$\phi(\mathbb{D}) \cap \mathbb{T} \ne \emptyset$ 
the operator $\phi(S^*)$ on $\ell^2 (\mathbb{N}_0)$ 
has a hypercyclic subspace.}
\bigskip

Note that the $\phi(z) = \lambda z$, $|\lambda| > 1$, 
does not satisfy this condition.

The examples in the Main Theorem may be interpreted as certain Toeplitz operator
on the Hardy space. 
The Hardy space $H^2 = H^2(\mathbb{D})$ is the space of all functions of the form 
$f(z) = \sum_{n\ge 0} c_n z^n$ with $\{c_n\} \in \ell^2(\mathbb{N}_0)$,
and thus is naturally identified with $\ell^2(\mathbb{N}_0)$. Recall
that for a function $\varphi\in L^\infty(\mathbb{T})$ the Toeplitz operator
$T_\varphi$ with the symbol $\varphi$ is defined as 
$T_\varphi f = P_+(\varphi f)$, where $P_+$ stands for the orthogonal projection from
$L^2(\mathbb{T})$ onto $H^2$. Then the backward shift on $S^*$ may be identified with 
the Toeplitz operator $T_{\overline z}$. It was shown in 
\cite{God-Sha} that any coanalytic Toeplitz operator $T_{\overline \varphi}$ 
(i.e., $\varphi$ is a bounded analytic function in $\mathbb{D}$) 
is hypercyclic whenever $\varphi(\mathbb{D})$
intersects $\mathbb{T}$. Our Main Theorem provides a class of 
coanalytic Toeplitz operators which have a hypercyclic subspace.

A general sufficient condition for the existence of a hypercyclic 
subspace was given by Gonzalez, Leon-Saavedra and Montes-Rodriguez in \cite{GLSMR}. To state it we need the following stronger 
version of hypercyclicity: 
\bigskip
\\
\textbf{Definition.} 
Operator $T$ acting on a separable Banach space $\mathcal{B}$ is 
said to be \textit{hereditarily hypercyclic} if there exists 
a sequence of non-negative integers $\{n_k\}$ 
such that for each subsequence $\{n_{k_i}\}$ there exists 
a vector $x$ such that the sequence $\{T^{n_{k_i}} x\}$ is 
dense in \bigskip $\mathcal{B}$.

We also need to recall the notion of the essential spectrum.
\bigskip
\\
\textbf{Definition.} Operator $U$ is called \textit{Fredholm} if $\mathrm{Ran}\ U$ is closed and has finite codimension and $\mathrm{Ker}\ U$ is finite-dimensional.
The \textit{essential spectrum} of the operator $T$ is defined as 
$$
 \sigma_e(T) = \{ \lambda: T - \lambda I\text{ is non-Fredholm}\}.
$$
\smallskip

\noindent
\textbf{Theorem (Gonzalez, Leon-Saavedra, Montes-Rodriguez, \cite[Theorem 3.2]{GLSMR}).} 
\textit{Let $T$ be a hereditary hypercyclic bounded linear operator on a separable 
Banach space $\mathcal{B}$. Let the essential spectrum of $T$ intersect the closed unit disc. Then there exists a hypercyclic subspace for the operator $T$.}
\bigskip

We intend to use this result in the proof of the Main Theorem.

Let us mention some other results on this topic.
Shkarin in \cite{Shk} proved that the differentiation operator 
on the standard Frechet space $\mathit{Hol} (\mathbb{C})$ has a hypercyclic subspace.
Quentin Menet in \cite[Corollary 5.5]{Men} 
generalized this result: he proved that for every non-constant polynomial $P$ the operator $P(D)$ has a hypercyclic subspace.
He also obtained some results concerning weighted shifts on $\ell^p$.

\bigskip

\bigskip


\section{Preliminaries on essential spectra of linear operators.}

The following lemma is well known. We give its proof for the convenience of the reader.
\bigskip
\\
\textbf{Lemma.} 
\textit{Essential spectrum of the operator $S^*$ is the unit circle.} 
\bigskip
\\
\textbf{Proof:} Let us consider three cases:
\smallskip

Case 1: $|\lambda| > 1$. Then the operator $S^* - \lambda I = 
-\lambda (I - \frac{1}{\lambda}S^*)$ is invertible and, thus, it is Fredholm.
\smallskip

Case 2: $|\lambda| < 1$. We have $S^* - \lambda I = S^* (I - \lambda S)$.
Since the operator $S^*$ is Fredholm (its kernel is one-dimensional, its 
image is the whole space $\ell^2$), and $I - \lambda S$ 
is invertible, their composition is also a Fredholm operator. 
\smallskip

Case 3: $|\lambda| = 1$. Then the operator $S^* - \lambda I$ is not Fredholm, because its image has infinite codimension.

Indeed, the pre-image of the sequence $(\lambda y_1, \lambda^2 y_2, \lambda^3 y_3,\lambda^4 y_4, \ldots) \in \ell^2$ 
is given by $(a, \lambda (y_1 + a), \lambda^2 (y_1 + y_2 + a), \ldots)$ and the 
equality $a = -\sum\limits_{i=1}^{+\infty} y_i$ is necessary for 
the inclusion of this sequence into $\ell^2$. 

Then the pre-image of the sequence 
\begin{equation}
\label{1}
\Big(1, \frac{1}{2}, \underbrace{0, \ldots, 0}_{\geq 2^2-1 \text{ times}}, \frac{1}{4}, \underbrace{0, \ldots, 0}_{\geq 2^4-1 \text{ times}}, \ldots, \frac{1}{2^n}, \underbrace{0, \ldots, 0}_{\geq 2^{2^n}-1 \text{ times}}, \ldots\Big), 
\end
{equation}
multiplied componentwise by $(\lambda, \lambda^2, \lambda^3, \ldots),$
is given by 
$$
\Big(-2, -1, \underbrace{-\frac{1}{2}, \ldots, -\frac{1}{2}}_{\geq 2^2 \text{ times}}, 
\underbrace{-\frac{1}{4}, \ldots, -\frac{1}{4}}_{\geq 2^4 \text{ times}}, \ldots, \underbrace{-\frac{1}{2^n}, 
\ldots, -\frac{1}{2^n}}_{\geq 2^{2^n} \text{ times}}, \ldots\Big),
$$ 
multiplied componentwise by $(1, \lambda, \lambda^2, \ldots)$, 
but such sequences do not belong to $\ell^2$. All 
sequences of the form \eqref{1}, as is easily seen, form an 
infinite-dimensional subspace in $\ell^2$. 
\qed
\bigskip

The following important theorem about the mapping of the essential spectra
can be found, e.g., in \cite[p. 107]{gold}.
\bigskip
\\
\textbf{Essential Spectrum Mapping Theorem.} 
\textit{For any bounded linear operator $T$ in a Hilbert space $H$ and for any
polynomial $P$ one has $\sigma_e(P(T)) = P(\sigma_e(T)).$}

\bigskip

\bigskip

\section{Proof of the Main Theorem}

In the proof of hereditary hypercyclicity of the operator  
$\phi(S^*)$ we will use the following well-known  criterion 
due to Godefroy and Shapiro \cite{God-Sha}
(for an explicit statement see, e.g., \cite[Theorem 3.1]{GEP}):
\bigskip
\\
\textbf{Theorem (Godefroy--Shapiro criterion).} 
\textit{Let $T$ be a bounded linear operator
in a separable Banach space. Suppose that the subspaces 
$$
\begin{aligned}
X_0 & = {\rm span}\{ x \in X: Tx = \lambda x \text{ for some }
\lambda \in \mathbb{C}, \ |\lambda| < 1 \}, \\
Y_0 & = {\rm span}\{ x \in X: Tx = \lambda 
x \text{ for some } \lambda \ in\mathbb{C}, \ |\lambda| > 1 \},
\end{aligned}
$$
are dense in $X$. Then $T$ is hereditarily hypercyclic.}
\bigskip
\\
\textbf{Proof of the Main Theorem}
We should verify two conditions of the theorem of Gonzalez, 
Leon-Saavedra and Montes-Rodriguez. 

Any function $\varphi$ from disc-algebra can be approximated 
uniformly in $\overline{\mathbb{D}}$ by a sequence 
of polynomials $P_n$. So $P_n(S^*)$ tends to $\varphi(S^*)$ 
in the operator norm. 

We need to show that $\sigma_e(\varphi(S^*))$ intersects the 
closed unit disc. 
Since $\phi(\mathbb{T}) \cap \mathbb{T} \ne \emptyset$,
there exist $\lambda, \, \mu \in \mathbb{T}$ such 
that $\varphi (\lambda) = \mu$. 
Then $\mu_n = P_n(\lambda)$ tend to $\mu$.
By the Essential Spectrum Mapping Theorem for any polynomial 
$P$ one has $\sigma_e(P(S^*)) = P(\sigma_e(S^*)) = P(\mathbb{T})$.
In particular, $\mu_n = P_n(\lambda) \in \sigma_e(P_n(S^*))$
for any $n$, and so $P_n(S^*) - \mu_n I$ is not Fredholm.

Since the set of Fredholm operators is open in the operator norm
(see, e.g., \cite[Theorem 4.3.11]{davies}), 
the set of non-Fredholm operators is closed, whence
the limit of $P_n(S^*) - \mu_nI$, which is equal to $\varphi(S^*) - \mu I$, 
is not Fredholm, and $\mu$ belongs to the essential spectrum of $\varphi(S^*)$.
The first condition of the theorem by Gonzalez, 
Leon-Saavedra and Montes-Rodriguez is verified. 

It is well known that the condition 
$\phi(\mathbb{D}) \cap \mathbb{T} \ne \emptyset$ 
implies that $\phi(S^*)$ satisfies the 
Godefroy--Shapiro criterion. Let us briefly recall this argument.

Recall that the point spectrum of $S^*$ equals
$\sigma_p(S^*) = \{\lambda: |\lambda| < 1\}$
and the eigenvector is given by 
$(1, \lambda, \lambda^2, \cdots) \in \ell^2(\mathbb{N}_0)$,
or, if we pass to the Hardy space $H^2(\mathbb{D})$
using the natural identification
of $H^2$ with $\ell^2(\mathbb{N}_0)$, by
$$
k_\lambda(z) = \frac{1}{1 - \overline{\lambda} z} = \sum_{n\ge 0}\lambda^n z^n.
$$ 
These are the Cauchy kernels, which are reproducing kernels of $H^2$.
Clearly, $k_\lambda$, $\lambda \in \mathbb{D}$, 
are also eigenvectors of $\varphi(S^*)$ which correspond to 
eigenvalues $\varphi(\lambda)$.

By the condition 
$\phi(\mathbb{D}) \cap \mathbb{T} \ne \emptyset$, we know
that $\varphi(\mathbb{D})$ is an open set 
which intersects both $\mathbb{D}$ and $\mathbb{C} 
\setminus \overline{\mathbb{D}}$.
Clearly, both of the sets $X_0 = \{k_\lambda, \lambda \in \mathbb{D}: 
|\varphi(\lambda)| > 1\}$ and 
$Y_0 = \{k_\lambda, \lambda \in \mathbb{D}: |\varphi(\lambda)| < 1\}$ are dense in $H^2$. 
Indeed, $f \in H^2$ is orthogonal to $k_\lambda$ if and only 
if $f(\lambda) = 0$ and both $\{\lambda \in \mathbb{D}: |\varphi(\lambda)| > 1\}$ 
and $\{\lambda \in \mathbb{D}: |\varphi(\lambda)| < 1\}$ are open sets.
Thus the conditions of the Godefroy--Shapiro criterion are 
satisfied and the hereditarily hypercyclicity of the operator $\phi(S^*)$ follows.

Thus, by the theorem
of Gonzalez, Leon-Saavedra and Montes-Rodriguez,  
the operator $\phi(S^*)$ has a hypercyclic subspace.
\qed
\bigskip
\\
{\bf Acknowledgements.}
The author is grateful to  Quentin Menet for helpful comments
on the first version of the paper.


\begin{thebibliography}{00.}

\bibitem{Bay-Mat}
\emph{F.~Bayart, E.~Matheron, \/}
Dynamics of Linear Operators.
Cambridge University Press, 2009.

\bibitem{Bou}
\emph{P.~S.~Bourdon, \/}
Invariant manifolds of hypercyclic vectors.
Proceedings of the American Mathematical Society, (3) \textbf{118} (1993), pp.~845--847.

\bibitem{davies}
\emph{E.~B.~Davies, \/}
Linear Operators and Their Spectra, 
Cambridge Studies in Advanced Advanced Mathematics, Vol. 106, 
Cambridge University Press, 2007.

\bibitem{God-Sha}
\emph{G.~Godefroy, J.~H.~Shapiro, \/}
Operators with dense, invariant, cyclic vector manifolds.
Journal of Functional Analysis, \textbf{98} (1991), pp.~229--269.

\bibitem{gold}
\emph{S.~Goldberg, \/} Unbounded Linear Operators, McGraw-Hill, New York, 1966.

\bibitem{GLSMR}
\emph{M.~Gonzalez, F.~Leon-Saavedra, A.~Montes-Rodriguez, \/} 
Semi-Fredholm Theory: Hypercyclic and supercyclic subspaces.
Proceedings of the London Mathematical Society, (3) \textbf{81} (2000), pp.~169--189.

\bibitem{Gro}
\emph{K.-G.~Grosse-Erdmann, \/} 
Universal families and hypercyclic operators. 
Bulletin of American Mathematical Society, (3) \textbf{36} (1999), pp.~345--381.

\bibitem{GEP}
\emph{K.-G.~Grosse-Erdmann, A.~Peris Manguillot, \/}
Linear Chaos, Springer, Berlin, 2011.

\bibitem{Men}
\emph{Q.~Menet, \/} 
Hypercyclic subspaces and weighted shifts.
arXiv:1208.4963v1 [math.FA], 24 Aug 2012.

\bibitem{MR}
\emph{A.~Montes-Rodriguez, \/}
Banach spaces of hypercyclic vectors.
Michigan Mathematical Journal, \textbf{43} (1996), pp.~419--436.

\bibitem{rol}
\emph{S.~Rolewicz, \/} On orbits of elements, Studia Math. {\bf 32} 
(1969), pp.~17--22.

\bibitem{Shk}
\emph{S.~Shkarin, \/}
On the set of hypercyclic vectors for the differentiation operator.
Israel Journal of Mathematics, \textbf{180} (2010), pp.~271--283.

\end{thebibliography}
\end{document}